\documentclass[11pt,reqno]{amsart}
\usepackage{amsfonts,amssymb,amsmath}
\usepackage{color}
\usepackage{hyperref}

\DeclareMathOperator{\bmax}{\mathbf{max}}
\DeclareMathOperator{\llangle}{\langle\!\langle\!{}}
\DeclareMathOperator{\rrangle}{{}\!\rangle\!\rangle}

\newtheorem{theorem}{Theorem}

\newtheorem{proposition}[theorem]{Proposition}

\newtheorem{remark}[theorem]{Remark}

\numberwithin{equation}{section} \numberwithin{theorem}{section}

\begin{document}

\title{Pattern Recognition on Oriented Matroids:
Symmetric Cycles in the Hypercube Graphs}

\author{Andrey O. Matveev}
\email{andrey.o.matveev@gmail.com}


\begin{abstract}
If $\mathfrak{V}$ is the vertex sequence of a symmetric $2t$--cycle in the hypercube graph with the vertices $\{1,-1\}^t$, then for any vertex $T$ of the graph
there exists a unique inclusion--minimal subset of $\mathfrak{V}$ such that $T$ is the sum of its elements. We present a simple combinatorial statistic on decompositions of vertices of the hypercube graphs with respect to symmetric cycles and describe their basic metric properties.
\end{abstract}

\maketitle

\pagestyle{myheadings}

\markboth{A.O.~MATVEEV}{PATTERN RECOGNITION ON ORIENTED MATROIDS}

\thispagestyle{empty}

\tableofcontents

\section{Introduction}

Let\hfill $\mathcal{H}:=(E_t,\{1,-1\}^t)$\hfill be\hfill the\hfill
oriented\hfill matroid,\hfill on\hfill the\hfill ground\newline set\hfill $E_t:=\{1,\ldots,t\}$,\hfill whose\hfill set\hfill of\hfill topes\hfill $\{1,-1\}^t$\hfill is\hfill the\hfill vertex\hfill set\hfill of \newline a~$t$--dimensional geometric hypercube in $\mathbb{R}^t$; we substitute the familiar sign components $+$ and $-$ of topes by the real numbers $1$ and $-1$ respectively. This oriented matroid is realizable as the {\em arrangement of coordinate hyperplanes\/} in $\mathbb{R}^t$, see~\cite[Example~4.1.4]{BLSWZ}.

The tope graph $\mathcal{T}(\mathcal{L}(\mathcal{H}))$
of the oriented matroid $\mathcal{H}$
is the {\em hypercube graph} with $2^t$ vertices, that is the {\em Hamming graph\/} $\boldsymbol{H}(t,2)$ closely related to the {\em Hamming association scheme\/} (also called the {\em $t$--cube}) $\mathbf{H}(t,2)$. See e.g.~\cite{Ba,B-I,BCN,CST,DGSh,DL,G,HIK,MS,O,St} on such graphs and association schemes.

The vertex set of the hypercube graph $\boldsymbol{H}(t,2)$ is the tope set of $\mathcal{H}$, that is\hfill the\hfill collection\hfill of\hfill all\hfill words\hfill (we\hfill regard\hfill them\hfill as\hfill row\hfill vectors\hfill from\hfill $\mathbb{R}^t$)\newline $T:=(T(1),\ldots,T(t))\in\{1,-1\}^t$; vertices $T'$ and $T''$ are adjacent in $\boldsymbol{H}(t,2)$ iff there is a unique component $e\in E_t$ such that $T'(e)=-T''(e)$.

Recall that the vertices of $\mathbf{H}(t,2)$ are in one--to--one correspondence with the elements of the Boolean lattice of subsets of the set $E_t$: it is convenient to regard the power set $\mathbf{2}^{E_t}$ of the set $E_t$ as the collection $\{T^-:\ T\in\{1,-1\}^t\}$
of the {\em negative parts\/} $T^-:=\{e\in E_t:\ T(e)=-1\}$ of topes of the oriented matroid $\mathcal{H}$.

Let $\boldsymbol{R}:=(R^0,R^1,\ldots,R^{2t-1},R^0)$ be a symmetric $2t$--cycle ({\em symmetric cycle}, for short) in $\boldsymbol{H}(t,2)$, that is, $R^{k+t}=-R^k$ for each $k$, $0\leq k\leq t-1$. The vertex sequence of any $(t-1)$--path in $\boldsymbol{R}$ is a {\em basis\/} of $\mathbb{R}^t$; indeed, the absolute value of the determinant of the matrix composed of those row vectors is $2^{t-1}$. This in particular means that the vertex sequence $\mathfrak{V}(\boldsymbol{R})$ $:=(R^0,R^1,\ldots,R^{2t-1})$ of the cycle $\boldsymbol{R}$  is a {\em maximal positive basis\/} of $\mathbb{R}^t$, see~\cite{AM1}.

For any vertex $T$ of $\boldsymbol{H}(t,2)$ there exists a unique inclusion--minimal subset~$\boldsymbol{Q}(T,\boldsymbol{R})\subset\mathfrak{V}(\boldsymbol{R})$ such that
\begin{equation*}
T=\sum_{Q\in\boldsymbol{Q}(T,\boldsymbol{R})}Q\; .
\end{equation*}
For all topes $T\in\{1,-1\}^t$, the linearly independent sets $\boldsymbol{Q}(T,\boldsymbol{R})\subset\mathbb{R}^t$ are of odd cardinality; they can be explicitly described as follows~\cite[Cor.~2.2]{AM1}:
\begin{equation*}
\boldsymbol{Q}(T,\boldsymbol{R}):={}_{-(T^-)}\!\left(\bmax^+\!\left({}_{-(T^-)}\mathfrak{V}(\boldsymbol{R})\right)\right)\; ,\ \ \ T\in\{1,-1\}^t\; ,
\end{equation*}
where ${}_{-(T^-)}\mathfrak{V}(\boldsymbol{R})$ is the sequence of topes obtained from the vertex sequence~$\mathfrak{V}(\boldsymbol{R})$ by {\sl reorientation\/} on the negative part $T^-$ of $T$; $\bmax^+(\boldsymbol{\cdot})$ is the subset of all topes from the resulting sequence that have {\sl inclusion--maximal positive parts}, and the outermost operation ${}_{-(T^-)}(\boldsymbol{\cdot})$ means the reverse {\sl reorientation\/} on the ground subset $T^-$.

If we consider the vertex set $\{1,-1\}^t$ of the hypercube graph $\boldsymbol{H}(t,2)$ as the disjoint union
$\dot\bigcup_{T\in\{1,-1\}^t}\{T\}$
of the singleton sets of its vertices, then
\begin{equation*}
\{1,-1\}^t=\dot\bigcup_{T\in\{1,-1\}^t}\Bigl\{\sum_{Q\in\boldsymbol{Q}(T,\boldsymbol{R})}Q\Bigr\}\; .
\end{equation*}
In other words, for any symmetric cycle $\pmb{R}$ in $\boldsymbol{H}(t,2)$ we have
\begin{equation*}
\{1,-1\}^t=\dot\bigcup_{S\subseteq E_t}\Bigl\{\sum_{Q\in{}_{-S}\!\left(\bmax^+\!\left({}_{-S}\mathfrak{V}(\boldsymbol{R})\right)\right)}Q\Bigr\}\; ;
\end{equation*}
note that if a pair $\{P',P''\}$ is the neighborhood of a word $P$ in the cycle~${}_{-S}\boldsymbol{R}$ then $P\in\bmax^+\!\left({}_{-S}\mathfrak{V}(\boldsymbol{R})\right)$ iff $|(P')^-|-1=|P^-|=|(P'')^-|-1$.

In Section~\ref{sec:exps} we find the cardinalities $|\boldsymbol{Q}(T,\boldsymbol{R})|$ of the decompositions of topes, in the context of arbitrary simple oriented matroids. In~Section~\ref{sec:mps} we list some metric relations for the sets $\boldsymbol{Q}(T,\boldsymbol{R})$. In Section~\ref{sec:scs} we describe a basic statistic associated with the vertices of hypercube graphs and their symmetric cycles.

\section{Decompositions of Topes with Respect to Symmetric Cycles in the Tope Graphs}

\label{sec:exps}

Let $\mathcal{M}:=(E_t,\mathcal{T})$ be a simple oriented matroid (it contains no loops, parallel or {\sl antiparallel} elements) with set of topes $\mathcal{T}$. Let us fix in the {\em tope\hfill graph}\hfill of\hfill $\mathcal{M}$\hfill a\hfill symmetric\hfill cycle\hfill $\boldsymbol{R}$\hfill with\hfill vertex\hfill sequence\hfill $\mathfrak{V}(\boldsymbol{R})$\newline $:=(R^0,R^1,\ldots,R^{2t-1})$. If $T\in\mathcal{T}$, then we define the row {\em distance vector~$\boldsymbol{z}_{T,\boldsymbol{R}}:=(z_{T,\boldsymbol{R}}(0),z_{T,\boldsymbol{R}}(1),\ldots,z_{T,\boldsymbol{R}}(2t-1))\in\ell^2(\mathbb{Z}_{2t})$ of the cycle $\boldsymbol{R}$ with respect to the tope\/} $T$ as follows:
\begin{equation*}
z_{T,\boldsymbol{R}}(k):=d(T,R^k)\; ,\ \ \ 0\leq k\leq 2t-1\; ,
\end{equation*}
where $d(T',T'')$ denotes the {\em graph distance\/} between topes $T'$ and $T''$, that is\hfill the\hfill {\em Hamming\hfill distance}\hfill between\hfill the\hfill words\hfill $T'$\hfill and\hfill $T''$.\hfill If\hfill we\hfill let\hfill $\langle T',T''\rangle$\newline $:=\sum_{e\in E_t}T'(e)T''(e)$ denote the standard scalar product on $\mathbb{R}^t$, then $\langle T',T''\rangle$ $=t-2d(T',T'')$ and, as a consequence, $d(T',T'')=\tfrac{1}{2}(t-\langle T',T''\rangle)$.

Let $\mathbf{I}$ and $\mathbf{C}$ denote the $2t\times 2t$ {\em identity matrix\/} and {\em basic circulant permutation matrix}, respectively, with the rows and columns indexed from $0$ to~$2t-1$. If $\boldsymbol{v}\in\ell^2(\mathbb{Z}_{2t})$, then we denote by $\hat{\boldsymbol{v}}:=(\hat{v}(0),\hat{v}(1),\ldots,\hat{v}(2t-1))$ the {\em discrete Fourier transform\/} (see~e.g.~\cite{F,W}) of the vector $\boldsymbol{v}$. If $\boldsymbol{w}\in\ell^2(\mathbb{Z}_{2t})$, then we let $\llangle\boldsymbol{v},\boldsymbol{w}\rrangle:=\sum_{k=0}^{2t-1}v(k)\overline{w(k)}$ denote the complex inner product on~$\ell^2(\mathbb{Z}_{2t})$, where $\overline{\phantom{\cdot}\!\!\boldsymbol{\cdot}\!\!\phantom{\cdot}}$
means complex conjugation.

For any tope $T\in\mathcal{T}$, we have
\begin{align}
\nonumber
|\boldsymbol{Q}(T,\boldsymbol{R})|&=t-\frac{1}{8}\boldsymbol{z}_{\!T,\boldsymbol{R}}
\cdot(2\mathbf{I}-\mathbf{C}^{-2}-\mathbf{C}^2)\cdot
\boldsymbol{z}_{\!T,\boldsymbol{R}}{}^{\!\top}\\
\label{eq:9}
&=
t-\frac{1}{4t}\sum_{k=0}^{2t-1} |\hat{z_{\!T,\boldsymbol{R}}}(k)|^2\cdot\sin^2\tfrac{\pi k}{t}
\; ,\\
\nonumber
|\boldsymbol{Q}(T,\boldsymbol{R})|&=
\frac{1}{8}\boldsymbol{z}_{\!T,\boldsymbol{R}}
\cdot(6\mathbf{I}-4\mathbf{C}^{-1}-4\mathbf{C}+\mathbf{C}^{-2}+\mathbf{C}^2)\cdot
\boldsymbol{z}_{\!T,\boldsymbol{R}}{}^{\!\top}\\
\label{eq:10}
&=
\frac{1}{4t}\sum_{k=0}^{2t-1} |\hat{z_{\!T,\boldsymbol{R}}}(k)|^2\cdot
(\cos^2\tfrac{\pi k}{t}-2\cos\tfrac{\pi k}{t}+1)\; ,\\
\nonumber
|\boldsymbol{Q}(T,\boldsymbol{R})|&=
\frac{t}{2}+\frac{1}{8}\boldsymbol{z}_{\!T,\boldsymbol{R}}
\cdot(2\mathbf{I}-2\mathbf{C}^{-1}-2\mathbf{C}+\mathbf{C}^{-2}+\mathbf{C}^2)\cdot
\boldsymbol{z}_{\!T,\boldsymbol{R}}{}^{\!\top}\\
\label{eq:11}
&=
\frac{t}{2}+
\frac{1}{4t}\sum_{k=0}^{2t-1} |\hat{z_{\!T,\boldsymbol{R}}}(k)|^2\cdot
(\cos^2\tfrac{\pi k}{t}-\cos\tfrac{\pi k}{t})
\intertext{and}
\nonumber
|\boldsymbol{Q}(T,\boldsymbol{R})|&=\frac{3t}{4}-\frac{1}{8}\boldsymbol{z}_{\!T,\boldsymbol{R}}
\cdot(\mathbf{C}^{-1}+\mathbf{C}-\mathbf{C}^{-2}-\mathbf{C}^2)\cdot
\boldsymbol{z}_{\!T,\boldsymbol{R}}{}^{\!\top}\\
\label{eq:12}
&=
\frac{3t}{4}+
\frac{1}{8t}\sum_{k=0}^{2t-1} |\hat{z_{\!T,\boldsymbol{R}}}(k)|^2\cdot
(2\cos^2\tfrac{\pi k}{t}-\cos\tfrac{\pi k}{t}-1)\; ,
\end{align}
see~\cite{AM2}.

Let $\boldsymbol{a}_{\!T,\boldsymbol{R}}\in\ell^2(\mathbb{Z}_{2t})$ denote the {\em autocorrelation\/} (see e.g.~\cite[Section~4.3]{Bu2nd}) of the distance vector $\boldsymbol{z}_{\!T,\boldsymbol{R}}$,
defined by
\begin{equation*}
a_{\!T,\boldsymbol{R}}(m):=\sum_{n=0}^{2t-1}z_{\!T,\boldsymbol{R}}(n)z_{\!T,\boldsymbol{R}}((n+m)\!\!\!\mod{2t})\; ,\ \ \  0\leq m\leq 2t-1\; ;
\end{equation*}
note\hfill that\hfill $a_{\!T,\boldsymbol{R}}(k)=a_{\!T,\boldsymbol{R}}(2t-k)$,\hfill $1\leq k\leq 2t-1$.\hfill Recall\hfill that\hfill $\hat{\boldsymbol{a}_{\!T,\boldsymbol{R}}}$\newline $=(|\hat{z_{\!T,\boldsymbol{R}}}(0)|^2,|\hat{z_{\!T,\boldsymbol{R}}}(1)|^2,\ldots,
|\hat{z_{\!T,\boldsymbol{R}}}(2t-1)|^2)$.

\begin{itemize}
\item[$\boldsymbol{\centerdot}$]
Let\hfill $\mathbf{b}:=(2,0,-1,0,\ldots,0,-1,0)$\hfill be\hfill the\hfill first\hfill row\hfill of\hfill the\newline matrix\hfill $2\mathbf{I}-\mathbf{C}^{-2}-\mathbf{C}^2$.\hfill  Eq.~(\ref{eq:9})\hfill is\hfill equivalent\hfill to\hfill the\newline
relation~\mbox{$t-|\boldsymbol{Q}(T,\boldsymbol{R})|=\tfrac{1}{16t}\llangle\hat{\boldsymbol{a}_{\!T,\boldsymbol{R}}},\hat{\mathbf{b}}\rrangle$,} and
Parseval's relation implies that
\begin{multline*}
t-|\boldsymbol{Q}(T,\boldsymbol{R})|=\tfrac{1}{8}\llangle\boldsymbol{a}_{\!T,\boldsymbol{R}},\mathbf{b}\rrangle\\=
\tfrac{1}{8}(2a_{\!T,\boldsymbol{R}}(0)-a_{\!T,\boldsymbol{R}}(2)-a_{\!T,\boldsymbol{R}}(2t-2))\\=
\tfrac{1}{8}(2a_{\!T,\boldsymbol{R}}(0)-2a_{\!T,\boldsymbol{R}}(2))
\; .
\end{multline*}

\item[$\boldsymbol{\centerdot}$]
Let\hfill $\mathbf{b}:=(6,-4,1,0,\ldots,0,1,-4)$\hfill be\hfill the\hfill first\hfill row\hfill of\hfill the\newline matrix $6\mathbf{I}-4\mathbf{C}^{-1}-4\mathbf{C}+\mathbf{C}^{-2}+\mathbf{C}^2$.
Eq.~(\ref{eq:10}) is equivalent to~$|\boldsymbol{Q}(T,\boldsymbol{R})|$ $=\tfrac{1}{16t}\llangle\hat{\boldsymbol{a}_{\!T,\boldsymbol{R}}},\hat{\mathbf{b}}\rrangle$;
Parseval's relation implies that
\begin{multline*}
|\boldsymbol{Q}(T,\boldsymbol{R})|=\tfrac{1}{8}\llangle\boldsymbol{a}_{\!T,\boldsymbol{R}},\mathbf{b}\rrangle\\=
\tfrac{1}{8}(6a_{\!T,\boldsymbol{R}}(0)-4a_{\!T,\boldsymbol{R}}(1)+a_{\!T,\boldsymbol{R}}(2)
+a_{\!T,\boldsymbol{R}}(2t-2)-4a_{\!T,\boldsymbol{R}}(2t-1))\\=
\tfrac{1}{8}(6a_{\!T,\boldsymbol{R}}(0)-8a_{\!T,\boldsymbol{R}}(1)+2a_{\!T,\boldsymbol{R}}(2))
\; .
\end{multline*}

\item[$\boldsymbol{\centerdot}$]
Let\hfill $\mathbf{b}:=(2,-2,1,0,\ldots,0,1,-2)$\hfill be\hfill the\hfill first\hfill row\hfill of\hfill the\newline matrix $2\mathbf{I}-2\mathbf{C}^{-1}
-2\mathbf{C}+\mathbf{C}^{-2}+\mathbf{C}^2$.
Eq.~(\ref{eq:11}) is equivalent to~$|\boldsymbol{Q}(T,\boldsymbol{R})|$ $-\tfrac{t}{2}=\tfrac{1}{16t}\llangle\hat{\boldsymbol{a}_{\!T,\boldsymbol{R}}},\hat{\mathbf{b}}\rrangle$, and we have
\begin{multline*}
|\boldsymbol{Q}(T,\boldsymbol{R})|-\tfrac{t}{2}=\tfrac{1}{8}\llangle\boldsymbol{a}_{\!T,\boldsymbol{R}},\mathbf{b}\rrangle\\=
\tfrac{1}{8}(2a_{\!T,\boldsymbol{R}}(0)-2a_{\!T,\boldsymbol{R}}(1)+a_{\!T,\boldsymbol{R}}(2)
+a_{\!T,\boldsymbol{R}}(2t-2)-2a_{\!T,\boldsymbol{R}}(2t-1))\\=
\tfrac{1}{8}(2a_{\!T,\boldsymbol{R}}(0)-4a_{\!T,\boldsymbol{R}}(1)+2a_{\!T,\boldsymbol{R}}(2))
\; .
\end{multline*}

\item[$\boldsymbol{\centerdot}$]
Let\hfill $\mathbf{b}:=(0,-1,1,0,\ldots,0,1,-1)$\hfill be\hfill the\hfill first\hfill row\hfill of\hfill the\newline matrix~$\mathbf{C}^{-1}+\mathbf{C}
-\mathbf{C}^{-2}-\mathbf{C}^2$.
Eq.~(\ref{eq:12}) is equivalent to $|\boldsymbol{Q}(T,\boldsymbol{R})|$ $-\tfrac{3t}{4}=\tfrac{1}{16t}\llangle\hat{\boldsymbol{a}_{\!T,\boldsymbol{R}}},\hat{\mathbf{b}}\rrangle$, and we have
\begin{multline*}
|\boldsymbol{Q}(T,\boldsymbol{R})|-\tfrac{3t}{4}=\tfrac{1}{8}\llangle\boldsymbol{a}_{\!T,\boldsymbol{R}},\mathbf{b}\rrangle\\=
\tfrac{1}{8}(-a_{\!T,\boldsymbol{R}}(1)+a_{\!T,\boldsymbol{R}}(2)
+a_{\!T,\boldsymbol{R}}(2t-2)-a_{\!T,\boldsymbol{R}}(2t-1))\\=
\tfrac{1}{8}(-2a_{\!T,\boldsymbol{R}}(1)+2a_{\!T,\boldsymbol{R}}(2))
\; .
\end{multline*}
\end{itemize}
We come to the conclusion:

\begin{proposition}
Let $\mathcal{M}$ be a simple oriented matroid, and $\mathfrak{V}(\boldsymbol{R})$
the vertex sequence of a
symmetric cycle $\boldsymbol{R}$ in the tope graph of $\mathcal{M}$.
If\/ $T$ is a tope of~$\mathcal{M}$,\hfill then\hfill for\hfill the\hfill inclusion--minimal\hfill subset\hfill $\boldsymbol{Q}(T,\boldsymbol{R})\subset\mathfrak{V}(\boldsymbol{R})$\hfill
such\newline that~$T=\sum_{Q\in\boldsymbol{Q}(T,\boldsymbol{R})}Q$, we have
\begin{align}
\label{eq:13}
|\boldsymbol{Q}(T,\boldsymbol{R})|&=t-\tfrac{1}{4}\bigl(a_{\!T,\boldsymbol{R}}(0)-a_{\!T,\boldsymbol{R}}(2)\bigr)\; ,\\
|\boldsymbol{Q}(T,\boldsymbol{R})|&=\tfrac{1}{4}\bigl(3a_{\!T,\boldsymbol{R}}(0)-4a_{\!T,\boldsymbol{R}}(1)+a_{\!T,\boldsymbol{R}}(2)\bigr)\; ,\\
|\boldsymbol{Q}(T,\boldsymbol{R})|&=\tfrac{t}{2}+\tfrac{1}{4}\bigl(a_{\!T,\boldsymbol{R}}(0)-2a_{\!T,\boldsymbol{R}}(1)+a_{\!T,\boldsymbol{R}}(2)\bigr)\\
\intertext{and}
|\boldsymbol{Q}(T,\boldsymbol{R})|&=\tfrac{3t}{4}+\tfrac{1}{4}\bigl(-a_{\!T,\boldsymbol{R}}(1)+a_{\!T,\boldsymbol{R}}(2)\bigr)\; ,
\end{align}
where $\boldsymbol{a}_{\!T,\boldsymbol{R}}$ is the {\em autocorrelation} of a distance vector of the cycle $\boldsymbol{R}$ with respect to the tope $T$.
\end{proposition}

\section{Basic Metric Properties of Decompositions}

\label{sec:mps}

In this section we consider a simple oriented matroid $\mathcal{M}:=(E_t,\mathcal{T})$ with distinguished symmetric cycle $\boldsymbol{R}$ in its tope graph. If $T\in\mathcal{T}$, then we have
\begin{equation*}
\sum_{Q\in\boldsymbol{Q}(T,\boldsymbol{R})}d(T,Q)=\sum_{Q\in\boldsymbol{Q}(T,\boldsymbol{R})}\frac{1}{2}(t-\langle T,Q\rangle)=\frac{1}{2}|\boldsymbol{Q}(T,\boldsymbol{R})|t-\frac{1}{2}\underbrace{\sum_{Q\in\boldsymbol{Q}(T,\boldsymbol{R})}\langle T,Q\rangle}_{\|T\|^2=t}\; .
\end{equation*}

\begin{remark} For any tope $T$ of $\mathcal{M}$ we have
\begin{align}
\label{eq:3}
\sum_{Q\in\boldsymbol{Q}(T,\boldsymbol{R})}d(T,Q)&=\frac{1}{2}\bigl(|\boldsymbol{Q}(T,\boldsymbol{R})|-1\bigr)t\\
\intertext{and}
\label{eq:2}
|\boldsymbol{Q}(T,\boldsymbol{R})|&=1+\frac{2}{t}\sum_{Q\in\boldsymbol{Q}(T,\boldsymbol{R})}d(T,Q)\; .
\end{align}
\end{remark}

Now suppose that $T\not\in\mathfrak{V}(\boldsymbol{R})$ and $\boldsymbol{Q}(T,\boldsymbol{R}):=\{Q^1,\ldots,Q^{|\boldsymbol{Q}(T,\boldsymbol{R})|}\}$.
Note that
\begin{multline*}
t=\|T\|^2=\langle T,T\rangle=\Bigl\langle\sum_{Q\in\boldsymbol{Q}(T,\boldsymbol{R})}Q,\sum_{Q\in\boldsymbol{Q}(T,\boldsymbol{R})}Q\Bigr\rangle\\=
|\boldsymbol{Q}(T,\boldsymbol{R})|t+2\sum_{1\leq i<j\leq|\boldsymbol{Q}(T,\boldsymbol{R})|}\bigl(t-2d(Q^i,Q^j)\bigr)\\=
|\boldsymbol{Q}(T,\boldsymbol{R})|t+2\tbinom{|\boldsymbol{Q}(T,\boldsymbol{R})|}{2}t-4\sum_{1\leq i<j\leq|\boldsymbol{Q}(T,\boldsymbol{R})|}d(Q^i,Q^j)\; .
\end{multline*}

\begin{remark} For a tope $T$ of $\mathcal{M}$, such that $T\not\in\mathfrak{V}(\boldsymbol{R})$, we have
\begin{align}
\sum_{1\leq i<j\leq|\boldsymbol{Q}(T,\boldsymbol{R})|}d(Q^i,Q^j)&=\frac{1}{4}\bigl(|\boldsymbol{Q}(T,\boldsymbol{R})|^2-1\bigr)t\; ,\\
|\boldsymbol{Q}(T,\boldsymbol{R})|&=\sqrt{1+4\tfrac{\sum_{i<j}d(Q^i,Q^j)}{t}}\; .\\
\intertext{Moreover,}
\sum_{1\leq i<j\leq|\boldsymbol{Q}(T,\boldsymbol{R})|}d(Q^i,Q^j)&=\frac{1}{2}(|\boldsymbol{Q}(T,\boldsymbol{R})|+1)\sum_{Q\in\boldsymbol{Q}(T,\boldsymbol{R})}d(T,Q)\; .
\end{align}
\end{remark}

\section{Symmetric Cycles in the Hypercube Graphs}

\label{sec:scs}

We now turn to an investigation of the combinatorial properties of the vertex set $\{1,-1\}^t$ of the hypercube graph $\boldsymbol{H}(t,2)$
and its symmetric cycles. Let $\boldsymbol{R}$ be a symmetric cycle, with its distance vectors $\boldsymbol{z}_{T\!,\boldsymbol{R}}$, $T\in\{1,-1\}^t$.

According to Eq.~(\ref{eq:13}), we have
\begin{align*}
\sum_{T\in\{1,-1\}^t}|\boldsymbol{Q}(T,\boldsymbol{R})|&=\sum_{T\in\{1,-1\}^t}\left(
t-\frac{1}{4}\sum_{n=0}^{2t-1}
z_{T\!,\boldsymbol{R}}(n)\bigl(z_{T\!,\boldsymbol{R}}(n)-z_{T\!,\boldsymbol{R}}((n+2)\!\!\!\!\mod{2t})\bigr)\right)\\
&=2^tt-\frac{1}{4}\sum_{n=0}^{2t-1}\sum_{T\in\{1,-1\}^t}
z_{T\!,\boldsymbol{R}}(n)\bigl(z_{T\!,\boldsymbol{R}}(n)-z_{T\!,\boldsymbol{R}}((n+2)\!\!\!\!\mod{2t})\bigr)\\
&=2^tt-\frac{1}{2}t\sum_{T\in\{1,-1\}^t}
z_{T\!,\boldsymbol{R}}(n)\bigl(z_{T\!,\boldsymbol{R}}(n)-z_{T\!,\boldsymbol{R}}((n+2)\!\!\!\!\mod{2t})\bigr)\\
&=t\Bigl(2^t-\frac{1}{2}\sum_{\substack{0\leq i,j\leq t}}\mathtt{p}_{ij}^2i(i-j)\Bigr)\; ,
\end{align*}
where for any vertices $X,Y\in\{1,-1\}^t$, such that $d(X,Y)=2$, the quantity~$|\{Z\in\{1,-1\}^t:\ d(Z,X)=i,\ d(Z,Y)=j\}|=:\mathtt{p}_{ij}^2$ is the same; this is an {\em intersection number\/} of the Hamming association scheme $\mathbf{H}(t,2)$, see~e.g.~\cite[\S21.3]{MS}.
Thus, we have
\begin{align*}
\sum_{T\in\{1,-1\}^t}|\boldsymbol{Q}(T,\boldsymbol{R})|&=t\Bigl(2^t-\frac{1}{2}\sum_{\substack{1\leq i\leq t;\\
j\in\{i-2, i+2\}:\\ 0\leq j\leq t
}}\binom{t-2}{\frac{i+j}{2}-1}\!\!\!\!\underbrace{\binom{2}{\frac{i-j}{2}+1}}_{1\text{\ when $j\in\{i-2, i+2\}$}}\!\!\!\!i(i-j)\Bigr)\\&=
t\Bigl(2^t-\frac{1}{2}\sum_{\substack{1\leq i\leq t;\\
j\in\{i-2, i+2\}:\\ 0\leq j\leq t
}}\binom{t-2}{\frac{i+j}{2}-1}i(i-j)\Bigr)\; .
\end{align*}
Since
\begin{equation*}
\mathfrak{s}(t):=\sum_{\substack{1\leq i\leq t;\\
j\in\{i-2, i+2\}:\ 0\leq j\leq t
}}\binom{t-2}{\frac{i+j}{2}-1}i(i-j)=2\mathfrak{s}(t-1)=2^t\; ,
\end{equation*}
we come to the following conclusion:
\begin{remark}
\label{th:2}
Let $\boldsymbol{R}$ be a symmetric cycle in the hypercube graph $\boldsymbol{H}(t,2)$.
\begin{itemize}
\item[\rm(i)] We have
\begin{equation}
\begin{split}
\sum_{T\in\{1,-1\}^t}|\boldsymbol{Q}(T,\boldsymbol{R})|&=\frac{t}{2}\sum_{\substack{1\leq i\leq t;\\
j\in\{i-2, i+2\}:\\ 0\leq j\leq t
}}\mathtt{p}_{ij}^2i(i-j)\\&=
\frac{t}{2}\sum_{\substack{1\leq i\leq t;\\
j\in\{i-2, i+2\}:\\ 0\leq j\leq t
}}\binom{t-2}{\frac{i+j}{2}-1}i(i-j)\; ,
\end{split}
\end{equation}
that is,
\begin{equation}
\label{eq:1}
\sum_{T\in\{1,-1\}^t}|\boldsymbol{Q}(T,\boldsymbol{R})|=2^{t-1}t
\end{equation}
---the number of edges in $\boldsymbol{H}(t,2)$.
\item[\rm(ii)]
Eqs.~{\rm(\ref{eq:1})} and~{\rm(\ref{eq:3})} yield
\begin{equation}
\label{eq:4}
\sum_{T\in\{1,-1\}^t}\sum_{Q\in\boldsymbol{Q}(T,\boldsymbol{R})}d(T,Q)=2^{t-2}(t-2)t\; .
\end{equation}
\end{itemize}
\end{remark}

If $j$ is an odd integer, $1\leq j\leq t$, then define
\begin{equation*}
c_j(t):=\bigl|\bigl\{T\in\{1,-1\}^t:\ |\boldsymbol{Q}(T,\boldsymbol{R})|=j\bigr\}\bigr|\; ;
\end{equation*}
thus, we have
\begin{equation*}
\sum_{\substack{1\leq j\leq t:\\ \text{\rm $j$ odd}}}c_j(t)=2^t\; .
\end{equation*}
Since
\begin{equation*}
\sum_{\substack{1\leq j\leq t:\\ \text{\rm $j$ odd}}}jc_j(t):=\sum_{T\in\{1,-1\}^t}|\boldsymbol{Q}(T,\boldsymbol{R})|=2^{t-1}t=\sum_{0\leq i\leq t} i\tbinom{t}{i}\; ,
\end{equation*}
by Remark~\ref{th:2}(i), the quantities $c_j(t)$ are
read off from Eq.~(\ref{eq:1}):
\begin{theorem} \label{th:3} Let $\boldsymbol{R}$ be a symmetric cycle, with vertex set $\mathfrak{V}(\boldsymbol{R})$, in the hypercube graph $\boldsymbol{H}(t,2)$. Consider, for the vertices $T\in\{1,-1\}^t$ of the graph~$\boldsymbol{H}(t,2)$, the inclusion--minimal subsets $\boldsymbol{Q}(T,\boldsymbol{R})\subset\mathfrak{V}(\boldsymbol{R})$ such that~$T=\sum_{Q\in\boldsymbol{Q}(T,\boldsymbol{R})}Q$.

For any odd integer $j$, $1\leq j\leq t$, we have
\begin{equation}
c_j(t):=\bigl|\bigl\{T\in\{1,-1\}^t:\ |\boldsymbol{Q}(T,\boldsymbol{R})|=j\bigr\}\bigr|=2\tbinom{t}{j}\; .
\end{equation}
In other words, the polynomial
\begin{equation*}
\boldsymbol{\gamma}_t(\mathrm{x}):=\sum_{\substack{1\leq j\leq t:\\ \text{\rm $j$ odd}}}
c_j(t)\mathrm{x}^j\; ,
\end{equation*}
in the variable $\mathrm{x}$, is
\begin{equation}
\boldsymbol{\gamma}_t(\mathrm{x})=2\sum_{\substack{1\leq j\leq t:\\ \text{\rm $j$ odd}}}\tbinom{t}{j}\mathrm{x}^j\; .
\end{equation}
\end{theorem}
We can define $\boldsymbol{\gamma}_1(\mathrm{x})$ as $2\mathrm{x}$. The polynomials $\boldsymbol{\gamma}_t(\mathrm{x})$, where $2\leq t\leq 10$, are collected in the following table:
\vspace{2mm}
{\small
\begin{center}
\begin{tabular}{ccc}
\hline
$t$ & $\quad$ &  $\boldsymbol{\gamma}_t(\mathrm{x}):=\sum_{T\in\{1,-1\}^t}\mathrm{x}^{|\boldsymbol{Q}(T,\boldsymbol{R})|}$
\\
\hline
\hline
$2$ & $\quad$ & $4\mathrm{x}$\\
$3$ & $\quad$ &  $6\mathrm{x}+2\mathrm{x}^3$\\
$4$ & $\quad$ &  $8\mathrm{x}+8\mathrm{x}^3$\\
$5$ & $\quad$ &  $10\mathrm{x}+20\mathrm{x}^3+2\mathrm{x}^5$\\
$6$ & $\quad$ &  $12\mathrm{x}+40\mathrm{x}^3+12\mathrm{x}^5$\\
$7$ & $\quad$ &  $14\mathrm{x}+70\mathrm{x}^3+42\mathrm{x}^5+2\mathrm{x}^7$\\
$8$ & $\quad$ &  $16\mathrm{x}+112\mathrm{x}^3+112\mathrm{x}^5+16\mathrm{x}^7$\\
$9$ & $\quad$ &  $18\mathrm{x}+168\mathrm{x}^3+252\mathrm{x}^5+72\mathrm{x}^7+2\mathrm{x}^9$\\
$10$ & $\quad$ &  $20\mathrm{x}+240\mathrm{x}^3+504\mathrm{x}^5+240\mathrm{x}^7+20\mathrm{x}^9$\\
\hline
\end{tabular}
\end{center}
}
\vspace{2mm}
In view of the simplicity of the statistic on the vertices and symmetric cycles of the hypercube graphs, given by the polynomials $\boldsymbol{\gamma}_t(\mathrm{x})$, there are many ways to describe the binomial--type combinatorial properties of the quantities\hfill $c_j(t)$.\hfill
For\hfill instance,\hfill
if\hfill $t$\hfill is\hfill even,\hfill then\hfill for\hfill any\hfill odd\hfill integer\hfill $j$,\newline $1\leq j< t$, we have $c_j(t)=c_{t-j}(t)$; if $t$ is odd, then for any odd integer $j$, $1\leq j\leq t$, we have $jc_j(t)=(1+t-j)c_{1+t-j}(t)$, and so on.

\section*{Appendix: Symmetric Cycles in the Hypercube Graph, and Equinumerous Decompositions of Vertices}

In this appendix we count the number of symmetric cycles $\boldsymbol{R}$ of the graph~$\boldsymbol{H}(t,2)$ that provide equinumerous decompositions~$\boldsymbol{Q}(T,\boldsymbol{R})$ of a vertex $T\in\{1,-1\}^t$.

It is easy to see that each of the $2^t$ vertices of the hypercube graph $\boldsymbol{H}(t,2)$ belongs to $\tfrac{1}{2}t!$ symmetric cycles. Indeed, the graph $\boldsymbol{H}(t,2)$ regarded as the~Hasse diagram of the {\em tope poset\/} of an oriented matroid (realizable as the arrangement of coordinate hyperplanes in $\mathbb{R}^t$) can be {\em based\/} at any tope~$B\in\{1,-1\}^t$, and such a poset, with its least element $B$ contained in $t!$ {\em maximal chains}, is a principal order ideal of a {\em binomial poset\/} whose {\em factorial function\/} is $n!$, see~\cite[Example~3.18.3b]{Stanley-EC}. It now suffices to note that any symmetric cycle of $\boldsymbol{H}(t,2)$ containing $B$ as its vertex is the union of two maximal chains of the tope poset based at $B$. Thus, there are~$\tfrac{1}{2t}\cdot 2^t\cdot\tfrac{1}{2}t!$ symmetric cycles in $\boldsymbol{H}(t,2)$:
\begin{equation*}
\#\{\boldsymbol{R}\colon \boldsymbol{R} \text{ symmetric cycle of $\boldsymbol{H}(t,2)$}\}=2^{t-2}(t-1)!\; .
\end{equation*}
See~\cite[\href{http://oeis.org/search?q=A002866}{A002866}]{OEIS} on the corresponding integer sequence.

Recall that for any symmetric cycle $\boldsymbol{R}$ of $\boldsymbol{H}(t,2)$ we have
\begin{equation*}
\sum_{T\in\{1,-1\}^t}|\boldsymbol{Q}(T,\boldsymbol{R})|=2^{t-1}t=2\sum_{\substack{1\leq j\leq t:\\ \text{$j$ odd}}}j\tbinom{t}{j}\; ,
\end{equation*}
see~Remark~\ref{th:2} and Theorem~\ref{th:3}. As a consequence, we have
\begin{equation*}
\sum_{\substack{\boldsymbol{R}:\\ \text{$\boldsymbol{R}$ symmetric cycle of $\boldsymbol{H}(t,2)$}}}\ \sum_{T\in\{1,-1\}^t}|\boldsymbol{Q}(T,\boldsymbol{R})|=2^{t-1}(t-1)!\sum_{\substack{1\leq j\leq t:\\ \text{$j$ odd}}}j\tbinom{t}{j}\; .
\end{equation*}
Since
\begin{equation*}
\begin{split}
\sum_{\substack{\boldsymbol{R}:\\ \text{$\boldsymbol{R}$ symmetric cycle of $\boldsymbol{H}(t,2)$}}}|\boldsymbol{Q}(T,\boldsymbol{R})|
&=\tfrac{1}{2^t}2^{t-1}(t-1)!\sum_{\substack{1\leq j\leq t:\\ \text{$j$ odd}}}j\tbinom{t}{j}\\&=
\tfrac{(t-1)!}{2}\sum_{\substack{1\leq j\leq t:\\ \text{$j$ odd}}}j\tbinom{t}{j}\; ,
\end{split}
\end{equation*}
for any vertex $T\in\{1,-1\}^t$ of the graph $\boldsymbol{H}(t,2)$,
we come to the following result:

\begin{proposition}
For any vertex $T\in\{1,-1\}^t$ of the hypercube graph~$\boldsymbol{H}(t,2)$, and for an {\em odd} integer $j$, $1\leq j\leq t$, we have
\begin{equation*}
\#\{\boldsymbol{R}\colon \boldsymbol{R} \text{ symmetric cycle of $\boldsymbol{H}(t,2)$},\ |\boldsymbol{Q}(T,\boldsymbol{R})|=j\}=
\tfrac{(t-1)!}{2}\tbinom{t}{j}\; .
\end{equation*}
\end{proposition}

\vspace{3mm}

\end{document}